%% file: modelFW.tex
\documentclass[12pt]{article}
 
\usepackage{myarticle}
\usepackage[colorinlistoftodos,prependcaption,textsize=tiny,backgroundcolor=black!5!white,bordercolor=red]{todonotes}

\usepackage{pgfplots}
\pgfplotsset{
  tick label style={font=\footnotesize},
  label style={font=\footnotesize},
  legend style={font=\footnotesize}
}

\usepackage{mycommands}
\usepackage{subcaption}
\usepackage{bm}

\newcommand{\eR}{(-\infty,+\infty]}        

\newcommand{\spC}{\mathscr C}
\renewcommand{\Jac}[1]{D\!#1}

\newcommand{\grow}{\omega}
\DeclareMathOperator{\diam}{diam}

\renewcommand{\proj}{\mathrm{proj}}
\newcommand{\matid}{I}

\newcommand{\mfun}[1]{f_{#1}}

\newcommand{\kp}{{k+1}}
\renewcommand{\k}{{k}}

\newcommand{\iter}[1]{_{#1}}
\newcommand{\piter}[1]{_{#1}}

\renewcommand{\rto}[1]{\overset{#1}{\to}}

\newcommand{\setsep}{\,\vert\,}

\graphicspath{{./figures/}}


\begin{document}

\thispagestyle{empty}
\begin{center}
\vspace*{0.03\paperheight}
{\Large\bf Model Function Based Conditional Gradient Method\\ \medskip with Armijo-like Line Search} \\
\bigskip
\bigskip
{\large Yura Malitsky$^\star$ and Peter Ochs$^\dagger$ \\ \medskip
{\small
$^\star$~Univeristy of G\"ottingen, G\"ottingen, Germany \\
$^\dagger$~Saarland University, Saarbr\"{u}cken, Germany \\
}
}
\end{center}
\bigskip

\begin{abstract}
The Conditional Gradient Method is generalized to a class of non-smooth non-convex optimization problems with many applications in machine learning. The proposed algorithm iterates by minimizing so-called model functions over the constraint set. Complemented with an Amijo line search procedure, we prove that subsequences converge to a stationary point. The abstract framework of model functions provides great flexibility for the design of concrete algorithms. As special cases, for example, we develop an algorithm for additive composite problems and an algorithm for non-linear composite problems which leads to a Gauss--Newton-type algorithm. Both instances are novel in non-smooth non-convex optimization and come with numerous applications in machine learning. Moreover, we obtain a hybrid version of Conditional Gradient and Proximal Minimization schemes for free, which combines advantages of both. Our algorithm is shown to perform favorably on a sparse non-linear robust regression problem and we discuss the flexibility of the proposed framework in several matrix factorization formulations.
\end{abstract}




\input{modelFW_intro}
\input{modelFW_rel}
\input{modelFW_alg}
\input{modelFW_modelex}
\newpage
\input{modelFW_num}

\section{Conclusion}

We have presented an algorithmic framework that generalizes the Conditional Gradient method from constrained convex or smooth minimization to a class of constrained non-smooth non-convex minimization problems. The algorithm is formulated with respect to sequential minimization of model functions over the constraint set, complemented with an Armijo line search procedure. Model functions are simple surrogates of the objective function that obey a certain approximation quality and capture first order information of the problem. We presented several examples of model functions, including examples for additive or non-linear composite problems, which demonstrates the gain in flexibility for solving problems in machine learning, computer vision, and statistics. The possibility to tailor model functions to the specific structure of the optimization problem at hand, allows for efficient minimization. We also devise a hybrid method that combines Conditional Gradient type update steps with proximal minimization steps, which is particularly interesting for matrix factorization problems. In a numerical experiment for robust non-linear regression, the algorithm performs favorably compared to proximal minimization based algorithms.

\subsection*{Acknowledgments}

This work was supported by the German Research Foundation (DFG) via grants SFB755-A4 and OC 150/1-1.

\input{modelFW_proofs}

{\small
\bibliographystyle{ieee}
\bibliography{ochs}
}

\end{document}

%% file: modelFW_intro.tex
\section{Introduction}

A prominent algorithm for applications in machine learning and statistics, such as matrix learning, recommender systems, clustering, etc., is the Conditional Gradient Method (aka Frank--Wolfe Method). Its success is based on a low per-iteration complexity in several applications. For example, in low rank approximation (e.g.\ matrix completion), the main computational cost per iteration is the minimization of a linear function over a nuclear norm (trace norm or Schatten 1-norm) constraint, which can be solved efficiently by approximating the singular vector associated with the largest singular value of the gradient that defines the linear function. In contrast, related proximal minimization algorithms require a full singular value decomposition, which is significantly more expensive. 

In this paper, we generalize the Conditional Gradient Method to non-smooth non-convex optimization problems and unify the convergence analysis for several algorithms. The classical convergence analysis relies on the Descent Lemma, in the case the objective $f$ has Lipschitz continuous gradient. The Descent Lemma states that
\[
  \abs{f(x) - f(\bar x) - \scal{\nabla f(\bar x)}{x-\bar x}} \leq \frac L2 \vnorm[2]{x-\bar x}^2 \quad \text{for all}\ x,\bar x\,.
\]
This inequality can also be interpreted as a measure for the linearization error of $f$ around $\bar x$, i.e., the approximation quality of $f$ by a linear function.  We emphasize the fact that such a measure for the approximation quality of $f$, rather than smoothness, is key for the convergence of the algorithm. We generalize the linear approximation to any \emphdef{model function} $\mfun{\bar x}$ that obeys a certain approximation quality 
\[
  \abs{f(x) - \mfun{\bar x}(x)} \leq \grow(\vnorm{x-\bar x})\,,
\]
measured by a growth function $\map{\grow}{\R_+}{\R_+}$ that controls the approximation error. Note that this inequality does not imply smoothness, even in the special case $\grow(t) = \frac{L}{2}t^2$. If $f=g+h$ with a smooth function $h$ and a non-smooth function $g$, we can define $\mfun{\bar x}(x) = g(x) + h(\bar x) + \scal{\nabla h(\bar x)}{x-\bar x}$ and observe that the approximation error is only due to the linearization of the smooth part $h$ of the objective, while $\mfun{\bar x}$ is non-smooth. There are many other situations of interest. We choose the properties of the growth function $\grow$ such that $\mfun{\bar x}$ mimics a first order oracle of $f$. The freedom to choose the model function depending on the problem structure at hand makes our approach a flexible and efficient way to solve structured non-smooth non-convex minimization problems.

In this model function framework, our generalized Conditional Gradient update step at $x\iter\k$ reads
\[
  \begin{split}
    y\iter\k \in&\ \argmin_{x\in C} \mfun{x\iter\k}(x) \\
    x\iter\kp =&\ \gamma\piter\k y\iter\k + (1-\gamma\piter\k) x\iter\k \,,
  \end{split}
\]
where $\gamma\piter\k\in[0,1]$ and $C$ is a compact and convex constraint set. For $\mfun{x\iter\k}$ being the linearization of $f$ around $x\iter\k$, this is exactly the Conditional Gradient Method. 

As for all methods, the efficiency depends on the cost to evaluate the oracle, which in our case is the minimization of $\mfun{x\iter\k}$ over $C$ and, for proximal minimization problems, the cost to solve subproblems of type
\[
  \min_{x\in \R^N} \mfun{x\iter\k}(x) + \frac{1}{2\tau}\vnorm{x-x\iter\k}^2\,,
\]
for some step size $\tau>0$.  The generalization achieved in this paper increases the modelling flexibility for practical applications by making them accessible with another (possibly much cheaper) oracle, or by combining the oracles to a hybrid Proximal--Conditional Gradient method. In particular, we show the favorable performance of our algorithm for a sparse non-linear robust regression problem and demonstrate the flexibility of the algorithm on several applications in matrix factorization.


%% file: modelFW_rel.tex
\section{Contributions and Related Work} 

The idea of model functions to unify and generalize algorithms has been used before in bundle methods \cite{Noll13, NPA08}, where only a lower bound on the approximation error with the model function is used, which is a different setup. In \cite{DIL16,OFB18}, the same class of model functions is considered as in our paper. In \cite{OFB18}, a Bregman proximal minimization framework is developed and convergence to a stationary point with an Armijo-like line search strategy is proved under weak assumptions on the Bregman distances. Their work can be seen as the proximal analogue to our framework. Recently, the model function framework has been extended to stochastic optimization \cite{DD18,DDM18}.

Both, \cite{OFB18} and our work, present an \emph{implementable algorithm of the model function framework}, which is motivated by the abstract consideration of (pure) sequential model minimization in \cite{DIL16}. The goal of \cite{DIL16} is to devise a measure for proximity to a stationary point, which can be used as a stopping criterion in non-smooth optimization. However, their convergence result depends on assumptions that are not automatically satisfied in practice. In \cite{OFB18}, model functions are complemented with additional structure (the Bregman proximity term) and an Armijo-like line search strategy. Once the model functions are selected, convergence of subsequences to a stationary point is guaranteed. We substitute the Bregman proximity by minimization of model functions over a compact set, and also obtain convergence of subsequences to a stationary point without additional assumptions. 

A special case of our framework yields the Conditional Gradient Method (aka Frank--Wolfe method \cite{FW56}) with Armijo line search. Convergence has been analyzed in \cite{Bertsekas99} for smooth constrained optimization and in \cite{RSPS16} for smooth stochastic problems. While, in convex optimization, convergence of the method is fairly well understood \cite{Bach15,Jaggi13,LJSP13,LJ15,SMF19,YFLC18,Nesterov18}, little is known in the non-smooth non-convex setting. To the best of our knowledge, our work is \emph{the first to generalize the Conditional Gradient minimization strategy to constrained non-smooth non-convex optimization with provable convergence (of subsequences) to a stationary point}. In this way, we contribute to the increase in modelling flexibility for problems in machine learning, computer vision, and statistics. In particular, we explore this \emph{flexibility in an example from non-linear robust regression and several formulations from matrix factorization}.

As specific instances of our algorithmic framework, we \emph{obtain new algorithms}. For example, we consider non-linear composite problems of type $\min_{x\in C} g(F(x))$ where $F$ is sufficiently smooth and $g$ is convex. Our iterative model function minimization over a convex constraint set yields an algorithm of Gauss-Newton type \cite{NW06}. Alternative strategies that use a proximal minimization strategy, which leads to Levenberg--Marquardt algorithm \cite{Levenberg44,Marquardt63} in a certain special case, is explored, for example, in \cite{LW16,DIL16,OFB18}. The problems that can be modelled in this form is immense \cite{LW16}. 
Using specific approximations of the objective by model functions, we also propose a \emph{hybrid Proximal--Conditional Gradient minimization scheme} that combines the advantages of both worlds. In the convex setting, such a hybrid method was used in \cite{ASS14}. However, their analysis was tailored to exactly this hybrid algorithm, whereas we obtain it from the model function framework for free and in the non-convex setting. 


%% file: modelFW_alg.tex
\section{Sequential Model Minimization with Line Search}
 
We consider optimization problems of the form
\begin{equation}\label{eq:objective}
  \min_{x\in C}\, f(x)
\end{equation}
with the following properties:
\begin{ASS} \label{ass:objective}
  \begin{itemize}
    \item[\ii1] $C$ is a non-empty compact convex set in $\R^N$;
    \item[\ii2] $\map{f}{\R^N}{\eR}$ is a proper lower semi-continuous (lsc) function that is bounded from below with $\dom f\subset C$.
    \end{itemize}
\end{ASS}
As motivated in the introduction, the proposed algorithm is based on iteratively minimizing model functions of the objective in \eqref{eq:objective} over the constraint set $C$. These model functions obey a certain approximation quality with respect to the objective function, which we measure in general using an (error) growth function:
\begin{DEF}[growth function] 
  A continuous function $\map{\grow}{\R_+}{\R_+}$ is called \emphdef{growth function} if it satisfies $\grow(0)=0$ and $\grow_+^\prime(0):=\lim_{t\dto 0} \grow(t)/t = 0$. 
\end{DEF}
The standard example of a  growth function is $\grow(t)=L\cdot t^r$ with $L>0$ and $r>1$. However, we may easily generate more examples using the concept of $\psi$-uniform continuity as in \cite{OFB18}, which generalizes Lipschitz and H\"older continuity.  
Note that $\grow(t) = o(t)$ if and only if $\grow$ is a growth function. 

In this paper, we consider model functions that satisfy the following assumption. 
\begin{ASS}[model assumption] \label{ass:model-assumption} 
  There exists a growth function $\map{\grow}{\R_+}{\R_+}$ such that for each $\bar x\in \R^N$, there exists a proper lsc convex function $\map{\mfun{\bar x}}{\R^N}{(-\infty, +\infty]}$ such that $\dom f=\dom\mfun{\bar x}$, called \emphdef{model function}, with the following property: 
  \[
       \abs{ f(x) - \mfun{\bar x}(x) } \leq \grow(\vnorm{x-\bar x})  \,, \qquad \forall x\in C\,.
   \]
\end{ASS}
For examples of model functions, we refer to
Section~\ref{sec:model-ex}. The Model
Assumption~\ref{ass:model-assumption} preserves up to the first order information of the objective function in the following sense (see Lemma~\ref{lem:differentials})
\begin{equation} \label{eq:FW-oracle}
  \mfun{\bar x}(\bar x)= f(\bar x) \quad\text{and}\quad \rpartial f(\bar x) = \partial \mfun{\bar x}(\bar x)  \,,
\end{equation}
where $\rpartial f$ denotes the Fr\'echet subdifferential \cite[Def. 8.3]{Rock98} of $f$ and $\partial f$ the (convex) subdifferential, which coincides with the Fr\'echet subdifferential for convex functions \cite[Prop. 8.12]{Rock98}.  The \emphdef{Fr\'echet subdifferential} is defined at a point $\bar x$, at which $f$ is finite, as $v\in\rpartial f(\bar x)$ if and only if $f(x) \geq f(\bar x) + \scal{v}{x-\bar x} + o(\vnorm{x-\bar x})$, and $\rpartial f(\bar x)=\emptyset$ for $\bar x\not\in\dom f$. \\

Minimizing model functions from Assumption~\ref{ass:model-assumption} provides a generic way to define algorithms with a first order oracle (possibly non-smooth). We seek to find a \emphdef{(Fr\'echet) stationary point} $\bar x$ of \eqref{eq:objective}, characterized by 
\[
  0 \in \rpartial f(\bar x)\,.
\]
In Algorithm~\ref{alg:modelFW}, the proposed algorithm is defined.

Key for measuring the progress of the algorithm is the \emphdef{model improvement}, which we define as
\begin{equation} \label{eq:def-model-improvement}
  \Delta(x,y) := \mfun{x}(x) - \mfun{x}(y) \quad \text{for all}\ x,y\in \R^N\,.
\end{equation}
Moreover, we show that it is a natural measure of stationarity.
\mybox{
  \begin{ALG}[Model Based Conditional Gradient Method with Line Search] \ \label{alg:modelFW}
    \vspace{-3ex}
    \begin{itemize}
      \item \key{Optimization Problem:} Problem \eqref{eq:objective}.
      \item \key{Initialization:} $x\iter0\in \R^N$ and set $\rho\in(0,1)$. 
      \item \key{Update $(\k\ge 0)$}: 
      \begin{itemize}
        \item Find $y\iter\k\in C$ such that the model improvement is positive, i.e.
        \begin{equation} \label{eq:FW-step}
          \Delta(x\iter\k,y\iter\k) = \mfun{x\iter\k}(x\iter\k) - \mfun{x\iter\k}(y\iter\k) > 0,
        \end{equation}
        and compute 
        \begin{equation} \label{eq:x-update}
          x\iter\kp = x\iter\k + \gamma\piter\k(y\iter\k - x\iter\k)
        \end{equation}
        with $\gamma\piter\k\in[0,1]$ determined by Algorithm~\ref{alg:armijo} such that the following holds:
        \begin{equation} \label{eq:armijo}
          \text{(Armijo line search)}\ 
          \gamma\piter\k \ \text{satisfies}\ f(x\iter\kp) \leq f(x\iter\k) - \rho\gamma\piter\k \Delta(x\iter\k,y\iter\k) \tag{ALS}
        \end{equation}
        \item If \eqref{eq:FW-step} cannot be satisfied (i.e., $\max_{y\in C}\Delta(x\iter\k,y)\leq 0$), then terminate the algorithm. 
      \end{itemize}
    \end{itemize}
  \end{ALG}
}
In order to obtain a ``stable'' algorithm, in the sense that objective values are non-increasing, the choice of $y\iter\k$ satisfying \eqref{eq:FW-step} is arbitrary. However, the proof that all limit points of the sequence generated by Algorithm~\ref{alg:modelFW} are stationary points requires an additional assumption. We must assert that the error in solving the model subproblem vanishes for $\k$ tending towards infinity.
\begin{ASS}[optimality of $y\iter\k$]\label{ass:y-opt}
  There exists $\seq[\k\in\N]{\eps\piter\k}$ with $\eps\piter\k\dto0$  such that 
  \[
      \mfun{x\iter\k}(y\iter\k) \leq \min_{x\in C} \mfun{x\iter\k}(x) + \eps\piter\k.
    \]
    For each $k\in \N$,  we denote by $\hat y\iter\k$ any element in $ \argmin_{x\in C} \mfun{x\iter\k}(x)$. 
\end{ASS}
  \begin{REM}
    One option to choose $y\iter\k$ in \eqref{eq:FW-step} is to set $y\iter\k = \hat y\iter\k\in \argmin_{x\in C} \mfun{x\iter\k}(x)$.  Observe that this is equivalent to $y_k \in \argmax_{y\in C}\Delta(x\iter\k,y)$. On the other hand, our framework is general enough to allow one solving $\min_{y\in C} \mfun{x\iter\k}(y)$ with errors as in Assumption~\ref{ass:y-opt}.
\end{REM}
The practical realization of the Armijo condition in \eqref{eq:armijo}, requires an algorithmic procedure. We propose the backtracking line search outlined in Algorithm~\ref{alg:armijo} as a subroutine for \eqref{eq:armijo}.
\mybox{
  \begin{ALG}[Armijo Line Search for Algorithm~\ref{alg:modelFW}] \ \label{alg:armijo}
    \begin{itemize}
      \item \key{Parameters:} Fix $\rho,\delta\in (0,1)$ and $\tilde \gamma\in (0,1]$.
      \item \key{Input:} $x\iter\k, y\iter\k \in C$ that satisfy \eqref{eq:FW-step}.
      \item \key{Line Search:} Find the smallest integer $j\geq0$ such that $\gamma\piter\k = \tilde\gamma \delta^j$ satisfies
        \eqref{eq:armijo}.
    \end{itemize}
  \end{ALG}
}
\subsection{Analysis of the Algorithm}

In the following sections, we discuss Algorithm~\ref{alg:modelFW}.

\subsubsection{Finite Termination of the Line Search Procedure}

We show that Algorithm~\ref{alg:modelFW} is well-defined, i.e., Algorithm~\ref{alg:armijo} terminates after a finite number of iterations. We verify that $y\iter\k-x\iter\k$ is a \emphdef{descent direction}, i.e., all sufficiently small choices of $\gamma\piter\k$ satisfy \eqref{eq:armijo}. Therefore, reducing $\gamma\piter\k$ according to the rule in Algorithm~\ref{alg:armijo}, it eventually enters a neighborhood of $0$ after finitely many steps. 
\begin{PROP} \label{prop:armijo-finite-term}
  Fix $\k\in \N$. There exists $\tilde \gamma\in (0,1]$ such that \eqref{eq:armijo} is satisfied for all $\gamma\piter\k\in (0,\tilde \gamma)$.
\end{PROP}
The proof is in Section~\ref{sec:proof:prop:armijo-finite-term}.

\subsubsection{Finite Termination of the Algorithm} 

In case, the algorithm terminates after a finite number of iterations, i.e., \eqref{eq:FW-step} cannot be satisfied for any $y\iter\k$, we have already found a stationary point.
\begin{PROP} \label{prop:Delta-stationarity}
  Let $\k\in \N$ be such that the model improvement is zero, i.e., $\max_{y\in C} \Delta(x\iter\k,y)=0$. Then, $x\iter\k$ is a stationary point of \eqref{eq:objective}. 
\end{PROP}
The proof is in Section~\ref{sec:proof:prop:Delta-stationarity}.\\

Proposition~\ref{prop:Delta-stationarity} identifies the model improvement $\Delta(x\iter\k,y\iter\k)$ as a suitable measure for stationarity. For smooth functions, this is an obvious fact, as the following example shows.
\begin{EX}
  If $f$ is sufficiently smooth, a suitable model function is $\mfun{x\iter\k}(x)=f(x\iter\k) + \scal{\nabla f(x\iter\k)}{x-x\iter\k}$, and the model improvement becomes
  \[
    \Delta(x\iter\k,y\iter\k) = \scal{\nabla f(x\iter\k)}{x\iter\k - y\iter\k} > 0 \,,
  \]
  which is the characterization of a descent direction $v=y\iter\k - x\iter\k$ in classical smooth optimization. If there is no $y\iter\k\in C$ along which the value of $\mfun{x\iter\k}$ can be reduced, then $\scal{\nabla f(x\iter\k)}{y-x\iter\k} \geq 0$ for all $y\in C$, which is the standard characterization of a stationary point $x\iter\k$ for constrained smooth optimization. 
\end{EX}

\subsubsection{Asymptotic Analysis}

In this section, we prove that all limit points of the sequence generated by Algorithm~\ref{alg:modelFW} are stationary points and that at least one such subsequence exists, which is stated in the following main convergence theorem.
\begin{THM}[convergence to a stationary point] \label{thm:convergence}
  Let Assumptions~\ref{ass:objective},~\ref{ass:model-assumption} and~\ref{ass:y-opt} be satisfied and let $\seq[\k\in\N]{x\iter\k}$ be a sequence that is generated by Algorithm~\ref{alg:modelFW}.  Then, every limit point of $\seq[\k\in\N]{x\iter\k}$ is a stationary point of \eqref{eq:objective} and $\seq[\k\in\N]{f(x\iter\k)}$ converges to the value of $f$ at the limit point. Moreover, there exists at least one converging subsequence of $\seq[\k\in\N]{x\iter\k}$. 
\end{THM}
The proof is in Section~\ref{sec:proof:thm:convergence}.
\begin{REM}
  Theorem~\ref{thm:convergence} guarantees to find a stationary point of the minimization problem in \eqref{eq:objective}. Note, that we do not intend to guarantee that a global minimizer of \eqref{eq:objective} is found or approximated. This would ask for too much considering the broadness of the class of non-smooth non-convex optimization problems that \eqref{eq:objective} deals with. In this general framework, convergence of subsequences to a stationary point is quite satisfying and is the objective of most first order optimization schemes in non-convex optimization. 
\end{REM}
\begin{REM}
  We can easily derive the following convergence rate from the Armijo line search condition \eqref{eq:armijo}:
  \[
      \min_{0\leq i \leq \k} \Delta(x\iter{i},y\iter{i}) \leq \frac{f(x\iter0) - \inf f}{\rho\sum_{i=0}^\k \gamma\piter\k}\,,\quad \forall \k\in\N\,.
  \]
  However, considering practice experiments, we observed that the convergence rate is too conservative and does not reflect the actual performance of our algorithm.
\end{REM}


%% file: modelFW_modelex.tex
\section{Examples of Model Functions} \label{sec:model-ex}

As the assumption of model functions (Assumption~\ref{ass:model-assumption}) is the same as in \cite{OFB18}, the same examples may be incorporated here. However, we consider minimization of model functions over the constraint set $C$ instead of (Bregman) proximal minimization. In order to make this work self contained, we mention their models (and some new ones) and discuss the algorithmic difference. For presentation, we focus on the case of maximal model improvement. Let $\Gamma_0$ denote the class of proper lsc convex functions and $\spC^{1,\psi}(C)$ be the class of smooth functions with $\psi$-uniformly continuous gradient relative to $C$, i.e., $f\in \spC^{1,\psi}$ if and only if
\[
  \vnorm{\nabla f(x) - \nabla f(y)} \leq \psi(\vnorm{x-y})\quad \text{for all}\  x,y\in C\,,
\]
for some continuous function $\map{\psi}{\R_+}{\R_+}$ with
$\psi(0)=0$. The Generalized Descent Lemma in \cite[Lem. 4]{OFB18}
shows that such a function obeys 
\begin{equation} \label{eq:gen-descent-lemma}
  \abs{f(x) - f(\bar x) - \scal{\nabla f(\bar x)}{x-\bar x}}\leq \grow(\vnorm{x-\bar x}) \quad \text{for all}\ x,\bar x\in C\,, 
\end{equation}
with a growth function $\grow(t):=\int_0^1 \frac{\phi(st)}{s}\, \mathit{d}s$ and $\map{\phi}{\R_+}{\R_+}$ given by $\phi(s) = s\psi(s)$. The most important example is $\psi(s)=cs^\alpha$ for some $c>0$, which is H\"older continuity for $\alpha\in(0,1]$ and Lipschitz continuity for $\alpha=1$. It results in $\grow(t) = \frac{c}{1+\alpha} t^{1+\alpha}$. Since the optimization problem in \eqref{eq:objective} is constrained to a compact set, local Lipschitz or H\"older continuity automatically become global (possibly with a different constant).

Note that in general the following examples account for non-smooth non-convex optimization problems.  
\begin{EX}[additive composite problems] \label{ex:fbs-model}
  Many problems in image processing, signal analysis, or statistics (including image deblurring, denoising, robust PCA, support vector machines, LASSO, etc.) can be cast in the form  
  \[
     \min_{x\in C} f(x)\,, \quad f:=g + h \quad \text{with}\  g\in \Gamma_0\ \text{and}\ h\in \spC^{1,\psi}(C) \,.
  \]
  A suitable model function for such problems is the following
  \[
     \mfun{\bar x}(x) = g(x) + h(\bar x) + \scal{\nabla h(\bar
     x)}{x-\bar x}, \quad \forall x \in C\,.
   \]
   Using \eqref{eq:gen-descent-lemma}, Assumption~\ref{ass:model-assumption} is clearly satisfied. 

  In the (Bregman) proximal minimization framework \cite{OFB18}, this choice requires to solve subproblems of the form 
  \[
    \min_{x\in C} g(x) + \scal{\nabla h(\bar x)}{x- \bar x} + D(x,\bar x) \,,
  \]
  which are known as (Bregman) Proximal Gradient Descent update steps (aka.\ Forward--Backward Splitting or Mirror Descent), where $D(x,\bar x)$ is a Bregman distance\footnote{The considered Bregman distances have the form $D(x,\bar x)=\phi(x) - \phi(\bar x) - \scal{\nabla \phi(\bar x)}{x-\bar x}$, if $\bar x\in \sint\dom \phi$, with a so-called Legendre function $\phi$, i.e., $\phi$ is essentially smooth and essentially strictly convex \cite[Sec. 26]{Rock70}, and $D(x,\bar x)=+\infty$ if $\bar x \not\in \sint\dom \phi$.}.  For $D(x,\bar x)=\frac{1}{2}\vnorm{x-\bar x}^2$ with $\tau>0$, the mapping that assigns to $\bar x$ the solution of this problem is known as the proximal gradient mapping $\prox\tau {g+\ind C}(\bar x - \tau\nabla h(\bar x))$ with respect to the function $g+\ind C$, where $\ind C$ is the indicator function of the set $C$.

  Instead, for our generalized Conditional Gradient type algorithm (Algorithm~\ref{alg:modelFW}) with maximal model improvement, the update step requires solving problems of type
  \[
     \min_{x\in C} g(x) + \scal{\nabla h(\bar x)}{x} \,.
  \]
  Key in selecting the ``better'' algorithm depends on the computational cost for solving the subproblems. 
\end{EX}
\begin{EX}[hybrid Proximal--Conditional Gradient minimization] \label{ex:prox-FW-model}
  Motivated by the comparison of proximal minimization and our Conditional Gradient type minimization in Example~\ref{ex:fbs-model}, the model function could be defined as 
  \[
    \mfun{\bar x}(x) = g(x) + h(\bar x) + \scal{\nabla h(\bar x)}{x-\bar x} + \frac{1}{2\tau}\vnorm{x-\bar x}^2\,,
  \]
  leading to a proximal subproblem over a constraint set $C$ in our Algorithm~\ref{alg:modelFW}. In this sense, our model function framework allows us to interpolate between proximal minimization algorithms and Conditional Gradient type algorithms. 
  
  We may also combine linearization and proximal linearization to devise a model function that yields a hybrid version of Conditional Gradient and proximal minimization. Consider the optimization problem in Example~\ref{ex:fbs-model} with $x=(x_1,x_2)\in \R^N$ and $C=C_1\times C_2$, where $g$ is additively separable, i.e., $g(x_1,x_2) = g_{1}(x_1) + g_{2}(x_2)$ for functions $g_1,g_2\in \Gamma_0$. Then, the following choice of model function
  \[
    \mfun{\bar x}(x) 
    = g_{1}(x_1) + g_{2}(x_2) + h(\bar x) + \scal{\nabla h(\bar x)}{x- \bar x} + \frac{1}{2\tau}\vnorm{x_1 - \bar x_1}^2 \,, \quad x=(x_1,x_2)\,,
  \]
  where $\nabla h(\bar x) = (\nabla_{x_1}h(\bar x), \nabla_{x_2}h(\bar x))$, leads to a proximal gradient step with respect to $x_1$ and a Conditional Gradient step with respect to $x_2$:
  \[
    \begin{split}
      \hat{y}_1 =&\ \prox{\tau}{g_{1}+\ind {C_1}}\Big( \bar x_1 - \tau \nabla_{x_1} h(\bar x) \Big) \\
      \hat{y}_2 \in&\ \argmin_{x_2\in C_2} g_{2}(x_2) + \scal{\nabla_{x_2} h(\bar x)}{x_2} 
    \end{split}
  \]
  where $\hat{y}=(\hat y_1, \hat y_2)$ yields the maximal model improvement in \eqref{eq:FW-step}.
\end{EX}
\begin{EX}[Newton-based Conditional Gradient] \label{eq:Newton-FW}
   Consider the problem in Example~\ref{ex:fbs-model} with higher regularity assumption, for example, suppose $h$ is at least twice continuously differentiable. In that case, a second order expansion in the model function is feasible 
   \[
      \mfun{\bar x}(x) = g(x) + h(\bar x) + \scal{\nabla h(\bar x)}{x-\bar x} + \frac 12 \scal{x-\bar x}{[\nabla^2 h(\bar x)]_+(x-\bar x)}\,,
   \]
   where $[\nabla^2 h(\bar x)]_+$ is the projection of the Hessian of $h$ at $\bar x$ onto the cone of positive semi-definite matrices. The convexity assumption of our model functions requires us to replace the Hessian matrix by a positive semi-definite approximation. However, in general, unlike proximal minimization methods, thanks to the compact constraint set, we do not need to enforce strong convexity of the subproblem, i.e., $[\nabla^2 h(\bar x)]_+$ need not be positive definite.  In the framework of \cite{OFB18} with $D(x,\bar x)=\frac{1}{2\tau}\vnorm{x-\bar x}^2$ and $g\equiv0$, this choice leads to damped (projected) Newton steps of the form
   \[
      \hat y = \proj_C \Big( \bar x - \tau (\matid + \tau [\nabla^2 h(\bar x)]_+)^{-1} \nabla h(\bar x) \Big)
   \]
   with identity matrix $\matid$. Our Algorithm~\ref{alg:modelFW} leads to a projected Newton step in subproblem \eqref{eq:FW-step} without damping
   \[
      \hat y = \proj_C\Big(\bar x - [\nabla^2 h(\bar x)]_+^{-1} \nabla h(\bar x)\Big) \,.
   \]
   Note the abuse of notation, since $[\nabla^2 h(\bar x)]_+$ might not be invertible. In that case, a constrained quadratic program needs to be solved to obtain a point $\hat y$ that yields the maximal model improvement. 
\end{EX}
\begin{EX}[Gauss--Newton] \label{ex:prox-linear}
  Consider minimization problems of the form
  \begin{equation} \label{eq:ex:prox-linear}
    \min_{x\in C} g(F(x)) \quad \text{with}\ g\in \Gamma_0(\R^M)\ \text{is Lipschitz and}\ F\in \spC^{1,\psi}(C,\R^M) \,.
  \end{equation}
  This class of problems includes non-linear inverse problems. We present a simple application of non-linear regression in Section~\ref{sec:robust-regression}. A suitable model function is the following:
  \[
    \mfun{\bar x}(x) = g(F(\bar x) + \Jac{F}(\bar x)(x-\bar x)) \,,
  \]
  which is motivated by the Gauss--Newton method \cite{NW06}. In the proximal minimization framework, it leads to the ProxLinear (or ProxDescent) algorithm \cite{LW16}, which can solve a broad class of problems. Often, the arising subproblems do not have closed form solution and numerical solvers are required. However, since the subproblems are convex, efficient minimization is possible. Due to the broad class of problems that is covered by \eqref{eq:ex:prox-linear}, in general, no simpler algorithms are currently known. There are essentially two ways of incorporating line search into the algorithm: \ii1 line search in direction of the solution of the subproblem \cite{OFB18} or \ii2 line search of the scaling of the proximity term to successively push the new iterate closer to the old iterate \cite{LW16,DIL16}. Where \ii1 requires to solve the subproblem once, \ii2 requires to solve the subproblem in each trial of a step size. 

  The line search strategy \ii1 is the same as in \eqref{eq:armijo}. As our subproblems do not involve the additional distance term, in contrast to proximal minimization subproblems, the search directions that we find are closer linked to the original problem, and hence we expect faster progress of our method. The robust regression problem in Section~\ref{sec:robust-regression} supports this intuition. 
\end{EX}
\begin{EX}
  The flexibility of our algorithm allows model functions to be tailored to specific problems. Suppose $g$ in \eqref{eq:ex:prox-linear} is additively separable, i.e., $g(y) = \sum_{i=1}^M g_i(y_i)$ and $g_i$ is convex and non-decreasing, e.g., the hinge loss $g_i(y_i) = \max(y_i, 0)$ that is used in support vector machines. Then, model functions with coordinate-wise higher order convex approximations of $F(x) = (F_1(x),\ldots,F_M(x))$ can be used to devise higher order convex model functions. 
\end{EX}


%% file: modelFW_num.tex
\section{Applications}

\subsection{Sparse Robust Non-linear Regression} \label{sec:robust-regression}

We consider a simple non-smooth non-convex sparse robust regression problem \cite{HRRS86} of the form
\[
  \min_{(a,b)\in C} \sum_{i=1}^M \vnorm[1]{F_i(a,b) - y_i} + \mu \vnorm[1]{a}\,, \quad F_i(a,b):= \sum_{j=1}^P a_j \exp(-b_j x_i) \,,
\]
similar to \cite{OFB18}, where the data $(x_i,y_i)\in \R^2$, $i=1,\ldots,M$, is a sequence of covariate-observation pairs and $C = [0,\overline a]^P \times [0,\overline b]^P$ for some $\overline a,\overline b>0$ and $P\in \N$. We assume that $y_i = F_i(a,b) + n_i$ where $n_i$ are iid errors drawn from a Laplacian distribution, which motivates the usage of the $\ell_1$-norm data fidelity term. Moreover, we assume that a large percentage of coefficients $a_j$ are zero, which is the reason for penalizing also the $\ell_1$-norm of the parameter vector $a\in \R^P$. By ``symmetry'' of $F_i$, the number of zero coordinates matters rather than the actual support.

We compare several algorithms with provable convergence of subsequences to a stationary point for solving the problem. The objective function falls into the class of problems of Example~\ref{ex:prox-linear}, for example, since $F$ has bounded Hessian on $C$, hence, its gradient is Lipschitz continuous on $C$. All algorithms are based on that choice of model functions. We write the linearization of the inner functions around $u\iter\k=(a\iter\k,b\iter\k)$ as follows: For all $i=1,\ldots,M$,
\[
  F_i(a,b) - y_i \approx \mathcal K_i u - y^\diamond_i\,,\quad \text{for all}\  u=(a,b)\in C\,, 
\]
where $\mathcal K_i = \Jac{F_i}(u\iter\k)$ and $y^\diamond_i:=y_i - F_i(u\iter\k) + \Jac{F_i}(u\iter\k) u\iter\k$.  Our Algorithm~\ref{alg:modelFW}, denoted \texttt{FW-CompLinLS}\footnote{Abbreviation for Frank--Wolfe Composite Linear splitting with line search.}, leads to subproblems of the form
\[
  \min_{u=(a,b)\in C} \sum_{i=1}^M \vnorm[1]{\mathcal K_i u - y_i^\diamond} + \mu \vnorm[1]{a} \,,
\]
the algorithm in \cite{OFB18}, denoted \texttt{ProxLinearLS}, and \cite{LW16}, denoted \texttt{ProxLinearBT}, require to solve subproblems of the form
\[
  \min_{u=(a,b)\in C} \sum_{i=1}^M \vnorm[1]{\mathcal K_i u - y_i^\diamond} + \mu \vnorm[1]{a} + \frac{1}{2\tau}\vnorm{u - u\iter\k}^2 \,.
\]
We solve the inner problem using the Primal--Dual Hybrid Gradient Algorithm with preconditioning \cite{PC11}, which allows for step sizes that are automatically computed based on $\mathcal K_i$. We use warm starting for all methods. Our algorithm \texttt{FW-CompLinLS} and \texttt{ProxLinearLS} solve the subproblem up to a certain accuracy and perform an Armijo-like line search in the direction of the approximate solution. The backtracking of \texttt{ProxLinearBT} is with respect to the parameter $\tau$ and involves solving the subproblem for each trial ``step size'' $\tau$ until a sufficient improvement of the objective value is observed. All methods perform line search for improving the objective value relative to the model improvement $\Delta(x\iter\k,y\iter\k)$. For details of the parameter setting, we refer to our code, which we provide with the paper. 

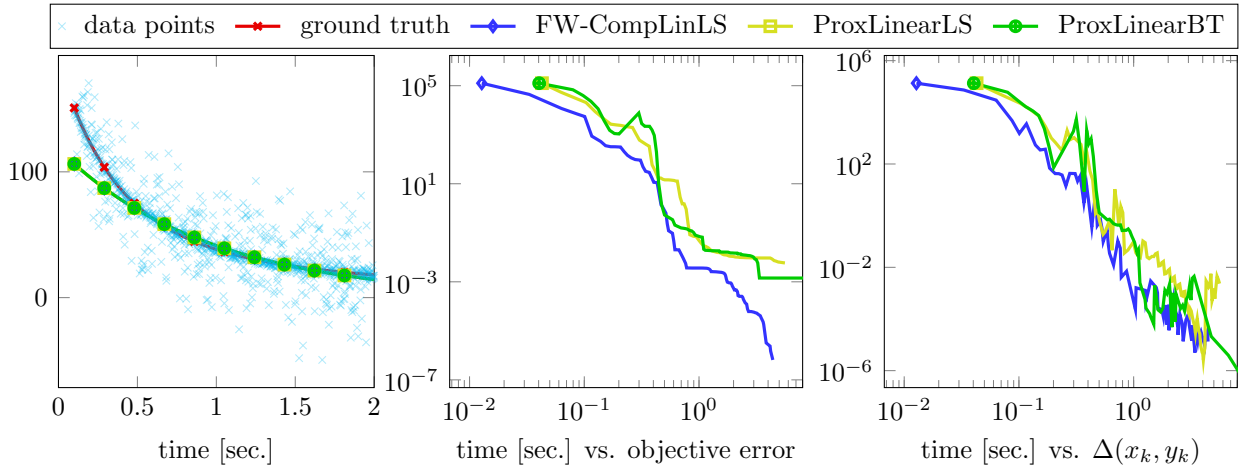
\begin{figure}[t]
  \begin{center}
    \begin{tikzpicture}
      \newcommand\markSize{2} 
      \begin{axis}[%
        width=0.35\linewidth,%
        height=6cm,%
        xmin=0,xmax=2,%
        xlabel={time [sec.]},%
        ylabel={},%
        legend columns=5,
        legend cell align=left,
        legend style={at={(1.855,1.02)},anchor=south,font=\footnotesize,column sep=5pt}
        ]
      \addplot[only marks,cyan!80!white,mark=x,opacity=0.3,mark size={\markSize},mark options={solid}] %
          table {figures/SignalAnalysis_data.dat};
          \addlegendentry{data points};
      \addplot[very thick,red!90!black,solid,mark=x,mark repeat={100},mark size={\markSize},mark options={solid}] 
          table {figures/SignalAnalysis_data_gt.dat};
          \addlegendentry{ground truth};
      \addplot[very thick,blue!80,solid,mark=diamond,mark repeat={100},mark size={\markSize},mark options={solid}] 
          table {figures/SignalAnalysis_regression_modelFW.dat};
          \addlegendentry{FW-CompLinLS};
      \addplot[very thick,yellow!80!green,solid,mark=square,mark repeat={100},mark size={\markSize},mark options={solid}] %
          table {figures/SignalAnalysis_regression_proxLinLS.dat};
          \addlegendentry{ProxLinearLS};
      \addplot[very thick,green!80!black,solid,mark=oplus,mark repeat={100},mark size={\markSize},mark options={solid}] 
          table {figures/SignalAnalysis_regression_proxLinBT.dat};
          \addlegendentry{ProxLinearBT};
      \end{axis}
      \begin{scope}[xshift=0.315\linewidth]
      \begin{axis}[%
        width=0.38\linewidth,%
        height=6cm,%
        xmin=0,xmax=8,%
        xlabel={time [sec.] vs. objective error},%
        ylabel={},%
        xmode=log, 
        ymode=log,
        ]
      \addplot[very thick,blue!80,solid,mark=diamond,mark repeat={100},mark size={\markSize},mark options={solid}] 
          table {figures/SignalAnalysis_obj_modelFW_time.dat};
      \addplot[very thick,yellow!80!green,solid,mark=square,mark repeat={100},mark size={\markSize},mark options={solid}] %
          table {figures/SignalAnalysis_obj_proxLinLS_time.dat};
      \addplot[very thick,green!80!black,solid,mark=oplus,mark repeat={100},mark size={\markSize},mark options={solid}] 
          table {figures/SignalAnalysis_obj_proxLinBT_time.dat};
      \end{axis}
      
      \end{scope}
          \begin{scope}[xshift=0.665\linewidth]
      \begin{axis}[%
        width=0.38\linewidth,%
        height=6cm,%
        xmin=0,xmax=8,%
        xlabel={time [sec.] vs. $\Delta(x\iter\k,y\iter\k)$},%
        ylabel={},%
        xmode=log, 
        ymode=log,
        ]
      \addplot[very thick,blue!80,solid,mark=diamond,mark repeat={100},mark size={\markSize},mark options={solid}] 
          table {figures/SignalAnalysis_Delta_modelFW_time.dat};
      \addplot[very thick,yellow!80!green,solid,mark=square,mark repeat={100},mark size={\markSize},mark options={solid}] %
          table {figures/SignalAnalysis_Delta_proxLinLS_time.dat};
      \addplot[very thick,green!80!black,solid,mark=oplus,mark repeat={100},mark size={\markSize},mark options={solid}] 
          table {figures/SignalAnalysis_Delta_proxLinBT_time.dat};
      \end{axis}
      \end{scope}
    \end{tikzpicture} 
  \end{center}
  \caption{\label{fig:rob-regr}Regression function and convergence plots for solving the robust regression problem in Section~\ref{sec:robust-regression}. All methods find the same regression function, as the left plot shows. The plot in the middle shows $f(x\iter\k)-\underline{f}$ where $\underline f$ is the smallest objective value found by any of the methods. The right plot shows the convergence of the model improvement, which is a measure for stationarity. The convergence is given with respect to actual computation time in seconds. Our method \texttt{FW-CompLinLS} outperforms the proximal line search \texttt{ProxLinearLS} and backtracking \texttt{ProxLinearBT} based methods.} 
\end{figure}

The data for the experiment is generated randomly with $P=100$, $M=1000$, $\mu=80$, $\overline{a}=20$, $\overline{b}=5$, and $80\%$ of coefficients $a_j$ are randomly set to $0$. Figure~\ref{fig:rob-regr} shows the data and the convergence of the objective value or the model improvement with respect to actual computation time.  

\subsection{Structured Matrix Factorization}

\newcommand{\XX}{\mathcal X}
\newcommand{\YY}{\mathcal Y}

Many applications in data analysis such as blind image deblurring \cite{KN06,CVR14}, clustering and principal component analysis \cite{DHS01,Murphy13}, source separation \cite{LS99,FBD09,CZPA09}, signal processing \cite{AEB06,SMF15},  or dictionary learning \cite{MBPS10,XLYZ17} can be formulated as structured matrix factorization problems. In this section, we demonstrate the flexible applicability of our algorithm to various formulations of matrix factorization problems. Most algorithms for solving such problems depend on alternating minimization techniques \cite{CZPA09,CVR14,SMF15}, sometimes with linearization \cite{BST14,PS16}. Algorithms are usually based on a proximal minimization oracle. In \cite{OFB18}, several formulations of matrix factorization are presented using Bregman proximal minimization steps. This approach has a great advantage for several constraint sets.

However, for example, proximal minimization of low rank constraints (e.g., constraints on the nuclear norm or 1-Schatten norm) require a full singular value decomposition (SVD), which can be expensive for large (or huge) scale data analysis problems \cite{CCS10}. In these settings, a Conditional Gradient minimization oracle is favorable. It requires to estimate the singular vector corresponding to the largest singular value only, which is computationally significantly cheaper than a full SVD. While this technique has been used frequently in (convex) low rank approximation schemes \cite{Jaggi13}, it has not been explored in detail for structured matrix factorization due to the non-convexity of the problem. We discuss several formulations of matrix factorization problems with focus on such low-rank constraints.  Due to the favorable properties of the generalized Conditional Gradient minimization oracle as described above, we believe that benchmarking is not required. 

We highlight the flexible applicability of our framework to non-convex
problems of the form
\[
    \min_{X,Y} \frac 12 \norm[F]{A - XY}^2 + g(X) \quad \st\ X\in \XX,\ Y \in \YY \,,
\]
where the goal is to represent a matrix $A$ as  a product $XY$ with
matrices $X\in \XX$ and $Y\in \YY$, where $\XX$ and $\YY$ are (convex and compact)
constraint sets that  encode  some problem specific characteristics and $g$ is a convex regularization function. We propose to use the additive composite splitting model from Example~\ref{ex:fbs-model}, i.e., we set $C= \XX\times \YY$,  $h(X,Y)=\frac 12 \norm[F]{A - XY}^2$ and solve the following subproblems:
\begin{equation} \label{eq:FW-update-MF}
    \min_{X,Y} g(X) + \scal{X}{(X\iter\k Y\iter\k-A)Y\iter\k^\top} + \scal{Y}{X\iter\k^\top(X\iter\k Y\iter\k - A)}_F\quad \st\ X\in \XX,\ Y \in \YY \,.
\end{equation}
Of course, the Frobenius norm in $h$ could be replaced by any smooth function, for example, the $\log$-student-t distribution $\sum_{i,j} \log(1+(A-XY)_{i,j}^2)$ for robust estimations. The linearization of $h$ makes the minimization separable, which allows us to discuss minimization steps with respect to $X$ and $Y$ independently. 

\paragraph{Examples for $\XX$.} In dictionary learning, $\XX$ describes the set of feasible atoms that may be used for reconstructing $A$. It is common to normalize the atoms, e.g.,
\[
  \XX_1 = \set{X \setsep \forall j\colon \sum_{i} X_{i,j}^2 \leq 1\,,\ \forall j>2\colon: \sum_{i} X_{i,j}=0}\,,
\]
which is a classical choice for dictionary learning \cite{XLYZ17}. For column $j=1$, the update step in \eqref{eq:FW-update-MF} is the projection of the $1$st column (of the gradient) onto the $\ell_2$-unit ball, and for $j>1$, by projecting the mean-subtracted $j$th column onto the $\ell_2$-unit ball. The choice
\[
  \XX_2 = \set{X \setsep \forall j\colon \sum_{i} X_{i,j} = 1\,,\ \forall i,j\colon X_{i,j}\geq 0} 
\]
enforces normalization and non-negativity, which is commonly used in non-negative matrix factorization (NMF) \cite{LS99}. The update step in \eqref{eq:FW-update-MF} sets column-wise a smallest coordinates to $1$ and all others to $0$. 

In \cite{OFB18}, a closed form update step with respect to $\XX_2$ is derived by a suitable choice of Bregman distance. Proximal minimization with respect to the Euclidean distance requires an algorithmic approach, though, which is also simple, as it is just a projection onto a unit simplex.

\paragraph{Examples for $\YY$ and $g$.} Sparsity is a favorable property for several matrix factorization problems. Conditional Gradient steps with respect to several norm constraints lead to simple updates \cite{Jaggi13,Bach13}. For example, for some $r>0$, set $\YY_1= \set{Y\setsep \norm[1]{Y}\leq r}$ to promote sparsity of the matrix $Y$. It can be used in dictionary learning \cite{XLYZ17} to express $A$ with only a few atoms of $X$, i.e., many entries of $Y$ shall be $0$. Analogously, convex relaxations of rank-$r$ constraints are commonly used, which can be modelled by $\YY_2 =\set{Y \setsep \norm[*]{Y} \leq r}$. The nuclear norm $\norm[*]{Y}$ of $Y$ enforces the columns of $A$ to be spanned by at most $r$ different linear subspaces, which is related to clustering problems. The Conditional Gradient subproblems with respect to both constraint sets $\YY_1$ and $\YY_2$ are simple \cite{Jaggi13}, where the second one requires the estimation of the extreme singular vector as mentioned above. 

However, on top of the constraint sets, we can use $g\not\equiv 0$, which may be used as penalty instead of a constraint, for example, penalizing the nuclear norm \cite{HDPDM12} or structured sparsity \cite{BJMO12}. The convex subproblem that arise in this context have been studied in convex optimization \cite{DHM12,HJN15,Nesterov18}. Also note that the solution of subproblems in \eqref{eq:FW-update-MF} with respect to $Y$ can be related to finding a subgradient in the subdifferential of the convex conjugate evaluated at the current gradient \cite{Bach15}.

\paragraph{Hybrid Proximal--Conditional Gradient minimization.} Finally, we discuss an alternative model function to \eqref{eq:FW-update-MF}, motivated by Example~\ref{ex:prox-FW-model}. We define the model function by linearization of the objective with respect to $Y$ and a convex quadratic approximation with respect to $X$. This choice leads to subproblems of the following form for our Algorithm~\ref{alg:modelFW}:
\[
    \min_{X,Y} g(X) + \scal{X}{(X\iter\k Y\iter\k-A)Y\iter\k^\top} + \scal{Y}{X\iter\k^\top(X\iter\k Y\iter\k - A)}_F + \frac{1}{2\tau}\vnorm[F]{Y - Y\iter\k}^2 \quad \st\ X\in \XX,\ Y \in \YY \,,
\]
for some $\tau>0$, leading to Conditional Gradient type problems with respect to $X$ and proximal minimization problems with respect to $Y$. The matrix factorization problem and the algorithm can be formulated to explore the advantages of both worlds. For example, a nuclear norm constraint with respect to $Y$ should be handled by a Conditional Gradient step and an additional group-sparsity penalty on $X$ can be efficiently handled by proximal minimization steps \cite{BJMO12}.


%% file: modelFW_proofs.tex
\appendix 

\section{Proofs}

\subsection{Proof of Proposition~\ref{prop:armijo-finite-term}} \label{sec:proof:prop:armijo-finite-term}

 For a fixed $k\in \N$, we abbreviate $\gamma=\gamma\piter\k$. Using Assumption~\ref{ass:model-assumption}, we have
  \[
    f(x\iter\kp) - f(x\iter\k) \leq \mfun{x\iter\k}(x\iter\kp) - \mfun{x\iter\k}(x\iter\k) + \grow(\vnorm{x\iter\kp-x\iter\k}) \,.
  \]
  From
  $\vnorm{x\iter\kp - x\iter\k} = \gamma\vnorm{y\iter\k -
  x\iter\k}$ and the definition of a growth
  function it follows that
  $\grow(\vnorm{x_{k+1}-x_k})=o(\gamma)$. The convexity of the model
  function $f_{x_k}$ gives us
  \[
      \mfun{x\iter\k}(x\iter\kp) - \mfun{x\iter\k}(x\iter\k) \leq \gamma (\mfun{x\iter\k}(y\iter\k) - \mfun{x\iter\k}(x\iter\k)) \,.
  \]
  Now, we argue by contradiction. Suppose that for any $\tilde \gamma >0$ there exists $\gamma\in (0,\tilde \gamma)$ such that \eqref{eq:armijo} does not hold, which yields the following calculation
  \[
    -\gamma \rho \Delta(x\iter\k, y\iter\k) < f(x\iter\kp) -
    f(x\iter\k) \leq \gamma (\mfun{x\iter\k}(y\iter\k) -
    \mfun{x\iter\k}(x\iter\k)) + o(\gamma) = -\gamma\Delta(x\iter\k, y\iter\k) + o(\gamma) \,.
  \]
  Dividing the inequality by $\gamma$, we obtain
  $0< (1-\rho) \Delta(x\iter\k, y\iter\k) <  o(\gamma)/\gamma$, which is a contradiction for sufficiently small $\tilde \gamma$.
\hfill\qedsymbol

\subsection{Proof of Proposition~\ref{prop:Delta-stationarity}}  \label{sec:proof:prop:Delta-stationarity}

The result is shown by Fermat's rule in the following lemma.
\begin{LEM} \label{lem:differentials}
  Let $\tilde x\in C$. Then,
  \[
    \rpartial f(\tilde x) = \partial \mfun{\tilde x}(\tilde x) \,,
  \]
  and 
  \[
    0 \in \partial \mfun{\tilde x}(\tilde x) \quad \Leftrightarrow\quad
    \Delta(\tilde x, x)  \leq 0 \ \forall x\in C\ .
  \]
\end{LEM}
\begin{proof}
  Let $v\in \rpartial f(\tilde x)$, then
  \[
    f(x) \geq f(\tilde x) + \scal{v}{x-\tilde x} + o(\vnorm{x-\tilde
    x}) \quad \forall x\in C
  \]
  and, this implies, by the model assumption
  \[
    \mfun{\tilde x}(x) + \grow(\vnorm{x-\tilde x}) \geq \mfun{\tilde x}(\tilde x) +\scal{v}{x-\tilde x} + o(\vnorm{x-\tilde x})\,, \quad \forall x\in C\,. 
  \]
  Since $\grow(t)=o(t)$, we conclude that
  \[
    \mfun{\tilde x}(x) \geq \mfun{\tilde x}(\tilde x)+  \scal{v}{x-\tilde x}  + o(\vnorm{x-\tilde x})\,, \quad \forall x\in C \,.
  \]
  Now, we fix a point $\bar x\in C$ and consider $x=\tilde x + \tau(\bar x - \tilde x)$ for $\tau \in (0,1]$. Then, by convexity of $C$ and the model function $\mfun{\tilde x}$, we obtain 
  \[
    \mfun{\tilde x}(\tilde x) + \tau (\mfun{\tilde x}(\bar x) - \mfun{\tilde x}(\tilde x)) \geq \mfun{\tilde x}(\tilde x) + \tau \scal{v}{\bar x- \tilde x} + o(\tau \vnorm{\bar x-\tilde x}).
  \]
  Subtracting $\mfun{\tilde x}(\tilde x)$, dividing by $\tau$, and considering $\tau \dto 0$, and, using the fact that this consideration was independent of the choice of $\bar x$, we conclude that $v\in\partial \mfun{\tilde x}(\tilde x)$. The converse direction follows easily. 

  The second part of the statement is Fermat's rule \cite[Thm 16.2]{BC11} for convex functions.
\end{proof}

\subsection{Proof of Theorem~\ref{thm:convergence}} \label{sec:proof:thm:convergence}

The result is proved in three steps.
\paragraph{Convergence of objective values.}
  The monotonicity and convergence of $\seq[\k\in\N]{f(x\iter\k)}$ follows directly from \eqref{eq:armijo} and the boundedness of $f$ from below.

\paragraph{Vanishing model improvement.}
  From \eqref{eq:armijo} and convergence of $\seq[\k\in\N]{f(x\iter\k)}$, we infer that $\gamma\iter\k \Delta(x\iter\k,y\iter\k)\to 0$, since
  \[
    0 \leq \rho \gamma\iter\k \Delta(x\iter\k,y\iter\k) \leq f(x\iter\k) - f(x\iter\kp) \to 0 \,.
  \]
  We deduce boundedness of $\seq[\k\in\N]{\Delta(x\iter\k,y\iter\k)}$ by
  \begin{multline*}
    0 
    \leq \Delta(x\iter\k,y\iter\k) 
    = \mfun{x\iter\k}(x\iter\k) - \mfun{x\iter\k}(y\iter\k) \leq f(x\iter\k) - \mfun{x\iter\k}(\hat y\iter\k) \\
    \leq f(x\iter0) - f(\hat y\iter\k) + \grow(\vnorm{\hat y\iter\k - x\iter\k}) 
    \leq f(x\iter0) - \inf_{x\in C} f(x) + \grow(\diam(C)) < +\infty.
  \end{multline*}
  Let $\Delta^*$ be an arbitrary limit point of $\seq[\k\in\N]{\Delta(x\iter\k,y\iter\k)}$, that is $\Delta(x\iter\k,y\iter\k)\to\Delta^*$ as $k\rto{K} \infty$ for some  $K\subset \N$, where $\k\rto{K}\infty$ abbreviates $\k\to\infty$ with $\k\in K$. 

  Suppose $\Delta^*>0$. Then $\gamma\piter\k\to 0$ as $ \k\rto{K}\infty$.  For sufficiently large $\k$, the line search procedure in Algorithm~\ref{alg:armijo} reduces $\gamma\piter\k/\delta$ to $\gamma\piter\k$, i.e., \eqref{eq:armijo} is violated before multiplying with $\delta$:
  \[
    -\frac{\gamma\piter\k}{\delta}\rho \Delta(x\iter\k, y\iter\k) < f(x\iter\k + \tfrac{\gamma\piter\k}{\delta}(y\iter\k-x\iter\k)) - f(x\iter\k) \,.
  \]
  Analogously to the proof of Proposition~\ref{prop:armijo-finite-term}, we conclude
  \[
    -\frac{\gamma\piter\k}{\delta}\rho \Delta(x\iter\k,y\iter\k) 
    < \frac{\gamma\piter\k}{\delta}(\mfun{x\iter\k}(y\iter\k) - \mfun{x\iter\k}(x\iter\k)) + o(\gamma\piter\k/\delta) 
    = -\frac{\gamma\piter\k}{\delta} \Delta(x\iter\k,y\iter\k) + o(\gamma\piter\k/\delta) \,.
  \]
  Dividing both sides by $\frac{\gamma\piter\k}{\delta}$ results in
  $(1-\rho) \Delta(x\iter\k, y\iter\k) <
  o(\gamma\piter\k)/\gamma\piter\k$ and considering $\gamma\piter\k\to
  0$ for $\k\rto{K}\infty$ yields a contradiction, since $\rho\in
  (0,1)$. Therefore $\Delta(x\iter\k,y\iter\k) \to 0$ for $\k\to\infty$.  

\paragraph{Convergence to a stationary point.}
  The following relation holds for all $x\in C$:
  \begin{equation} \label{eq:proof-prop:conv-pointwise-model-conv}
    \begin{split}
      \Delta(x\iter\k,y\iter\k) 
      =&\  \Delta(x\iter\k,\hat{y}\iter\k) + \mfun{x\iter\k}(\hat{y}) -  \mfun{x\iter\k}(y_k)  \\
      \geq&\  \mfun{x_k}(x_k) - \mfun{x_k}(x) - \eps\piter\k \\
      \geq&\ f(x_k) - f(x) - \grow(\vnorm{x_k-x}) - \eps\piter\k \,,
    \end{split}
  \end{equation}
  where the first inequality follows from Assumption~\ref{ass:y-opt} and the second from   Assumption~\ref{ass:model-assumption}.  Taking the limit $\k\rto{K}\infty$ on both sides, using $\Delta(x\iter\k,y\iter\k) \to 0$ for $\k\to\infty$, lower semi-continuity of $f$ and continuity of $\grow$, we arrive at
  \[
    f(x) \geq f(\tilde x) -\grow(\vnorm{\tilde x - x})\,, \quad \forall x\in C\,,
  \]
  where $\tilde x\in C$ due to compactness of $C$. As $\tilde x \in C$ and $\grow(t)=o(t)$, we deduce that 
  \[
    \liminf_{\substack{x \to \tilde x\\x\neq \tilde x}}\frac{f(x) - f(\tilde x) - \scal{0}{x-\tilde x}}{\vnorm{x-\tilde x}}\geq 0\,.
\]
which by definition means that $0\in \rpartial f(\tilde x)$.

Moreover, using $x=\tilde x$ in \eqref{eq:proof-prop:conv-pointwise-model-conv}, taking the limit $\k\rto{K}\infty$ and using lower semi-continuity of $f$, we deduce
\[
  f(\tilde x) \geq \limsup_{\k\rto{K}\infty} f(x\iter\k) \geq \liminf_{k\rto{K}\infty} f(x\iter\k) \geq f(\tilde x)\,,
\]
hence $f(x\iter\k) \to f(\tilde x)$ as $\k\rto{K}\infty$. By convergence of $\seq[\k\in\N]{f(x\iter\k)}$, we also have $f(x\iter\k)\to f(\tilde x)$ for $\k\to\infty$.
\hfill\qedsymbol
